\pdfoutput=1

\documentclass[letterpaper, 12pt, reqno]{amsart}
\usepackage{amsmath, amsfonts, amssymb, amsthm, array, stmaryrd, graphicx, hyperref, mathrsfs, eucal, caption, subcaption, soul, cancel, ulem, float}

\usepackage[alphabetic]{amsrefs}
\AtBeginDocument{%
	\def\MR#1{}
}

\usepackage{mathptmx}
\DeclareSymbolFont{Symbols}{OMS}{zplm}{m}{n}
\DeclareMathSymbol{\infty}{\mathord}{Symbols}{"31}

\usepackage[hhmmss]{datetime}

\allowdisplaybreaks

\usepackage[letterpaper, total={6.5in, 8.5in}, left=1in, top=1.25in]{geometry}
\newcount\Comments  
\usepackage{color}
\definecolor{burntorange}{RGB}{204, 85, 0}
\definecolor{oregongreen}{RGB}{0, 112, 48}
\definecolor{oaklandgold}{RGB}{181, 154, 87}
\definecolor{sydneyblue}{RGB}{1, 72, 164}

\newcommand{\kibitz}[2]{\ifnum\Comments=1\textcolor{#1}{#2}\fi}

\theoremstyle{plain}

\theoremstyle{definition}

\theoremstyle{remark}

\numberwithin{equation}{section}

\newcommand{\p}{\partial}

\newcommand{\bn}{\begin{enumerate}}
	\newcommand{\en}{\end{enumerate}}
\newcommand{\bi}{\begin{itemize}}
	\newcommand{\ei}{\end{itemize}}
\newcommand{\bqq}{\begin{eqnarray*}}
	\newcommand{\eqq}{\end{eqnarray*}}
\newcommand{\balg}{\begin{align*}}
\newcommand{\ealg}{\end{align*}}

\begin{document}
\title[Numerical stability analysis of noncompact Type-II MCF solutions: II]
{A numerical stability analysis of mean curvature flow of noncompact hypersurfaces with Type-II curvature blowup: II}

\author{David Garfinkle}
\address{Department of Physics, Oakland University, Rochester, MI 48309, USA}
\email{garfinkl@oakland.edu}

\author{James Isenberg}
\address{Department of Mathematics, University of Oregon, Eugene, OR 97403, USA}
\email{isenberg@uoregon.edu}

\author{Dan Knopf}
\address{Department of Mathematics, The University of Texas, Austin, TX 78712, USA}
\email{danknopf@math.utexas.edu}

\author{Haotian Wu}
\address{School of Mathematics and Statistics, The University of Sydney, NSW 2006, Australia}
\email{haotian.wu@sydney.edu.au}





\begin{abstract}

In previous work~\cite{GIKW}, we have presented evidence from numerical simulations that the Type-II singularities of mean curvature flow (MCF) of rotationally-symmetric, complete, noncompact embedded hypersurfaces constructed in \cite{IW, IWZ1} are stable. More precisely, it is shown in that paper that for small rotationally-symmetric perturbations of initial 
embeddings near the ``tip'', numerical simulations of MCF of such initial embeddings develop the same Type-II singularities with the same ``bowl soliton'' blowup behaviors in a neighborhood of the singularity. It is also shown in that work that for small
rotationally-symmetric perturbations of the initial embeddings that are sufficiently far away from the tip, MCF develops Type-I ``neckpinch'' singularities.

In this work, we again use numerical simulations to show that MCF subject to initial perturbations that are \emph{not} rotationally symmetric behaves asymptotically like it does for rotationally-symmetric perturbations. In particular, if we impose sinusoidal angular dependence on the initial embeddings, we find that for perturbations near the tip, evolutions by MCF asymptotically lose their angular dependence --- becoming round --- and develop Type-II bowl soliton singularities. As well, if we impose sinusoidal angular dependence on the initial embeddings for perturbations sufficiently far from the tip, the angular dependence again disappears as Type-I neckpinch singularities develop. The numerical analysis carried out in this work is an adaptation of the ``overlap'' method introduced in~\cite{GIKW} and permits angular dependence.

\end{abstract}


\maketitle

\tableofcontents



\section{Introduction}\label{intro}

In previous work \cite{GIKW}, we have used numerical simulations to show that Type-II singularity formations observed in mean curvature flow (MCF) of certain noncompact rotationally-symmetric embedded hypersurfaces~\cite{IW, IWZ1, IWZ2} are stable for small perturbations near the tip of each initial embedding, as long as rotational symmetry is retained. The question then arises whether this behavior is stable for small perturbations that are \emph{not} rotationally symmetric. In this paper, we present numerical simulations which indicate that Type-II singularity behaviors are also stable for such perturbations.

\begin{figure}[H]
	\includegraphics[width=0.65\textwidth]{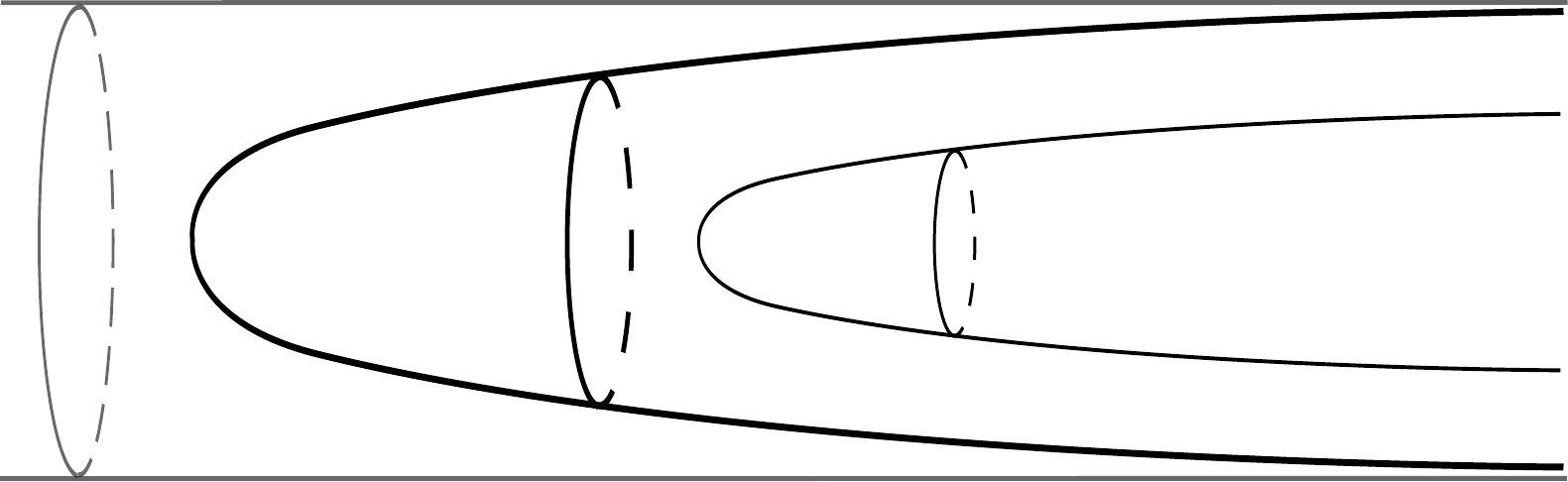}
	\caption{A MCF solution forming a degenerate neckpinch at spatial infinity with the curvature
	at the tip (on the left) blowing up at a Type-II rate \cite{IW,IWZ1}.}				\label{fig1}
\end{figure}

One says that MCF of an embedded hypersurface develops a Type-II singularity if the supremum of the product of the time to the singularity and the norm of the second fundamental form becomes infinite. If, alternatively, for MCF of a given embedded hypersurface, this product remains finite for the lifetime of the flow, then the singularity is defined to be Type-I. While Type-I singularities are expected to predominate for MCF (the same is expected to be true for Ricci flow), the work in \cite{IW, IWZ1, IWZ2} shows that Type-II singularities occur for MCF solutions originating from an open set of initial embeddings within the class of rotationally-symmetric noncompact hypersurfaces satisfying certain conditions. In these solutions, the maximum of the norm of the second fundamental form occurs at the tip (the left-most point on the  hypersurface depicted in Figure~\ref{fig1}), with the rate of curvature blowup at the tip consistent with the definition of a Type-II singularity. In \cite{GIKW}, we consider numerical simulations of MCF originating from rotationally-symmetric embedded hypersurfaces of this type having ``dimple'' perturbations of two classes.

\vspace{12pt}

\begin{figure}[H]
	\includegraphics[width=0.65\textwidth]{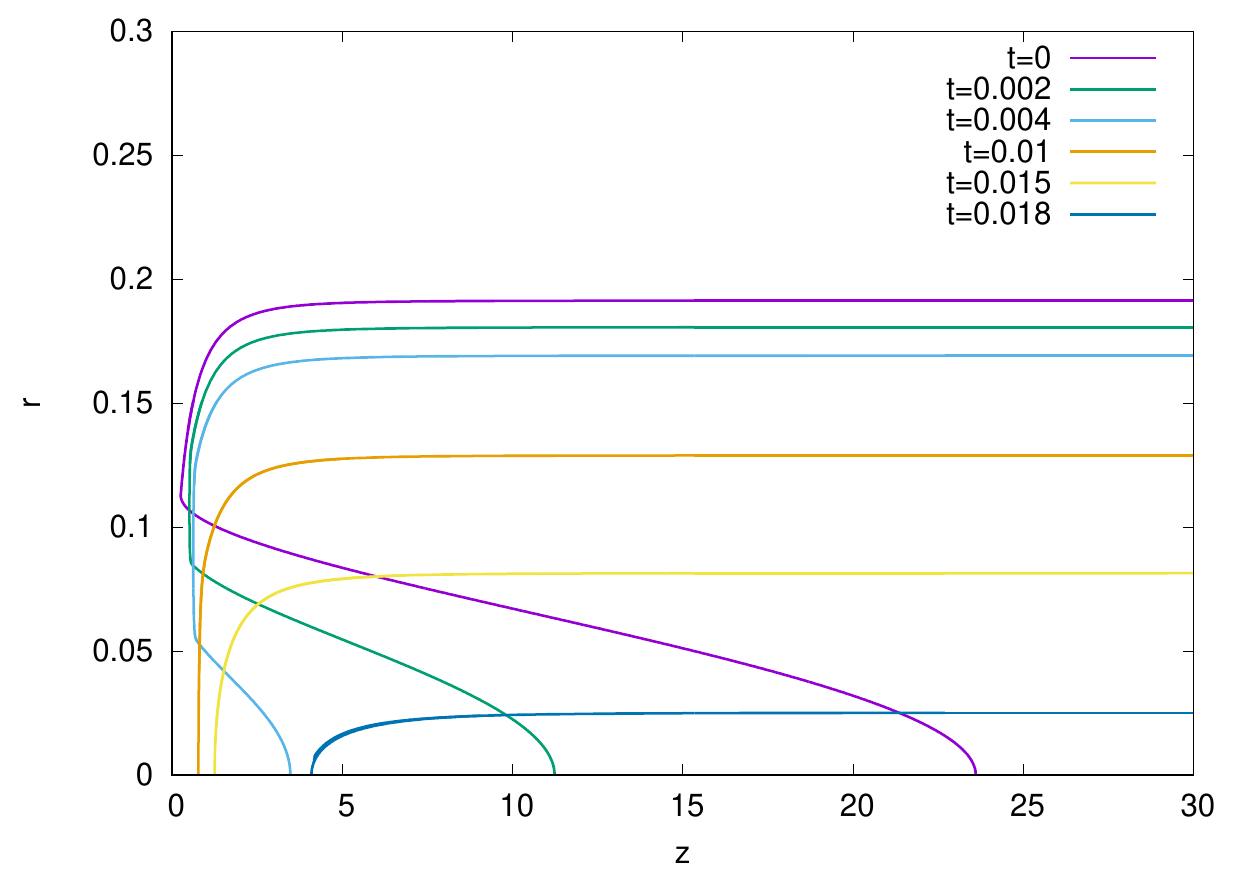}
	\caption{Numerical simulation of a MCF solution in the Near Class.
	Rotating each colored curve around the $z$-axis generates the hypersurface at the corresponding time.
	(The coordinates $z$ and $r$ are defined on page~\pageref{page:define-Coordinates} below.)
	As the Type-II singularity develops, the dimple disappears.}		\label{fig2}
\end{figure}

The ``Near Class'' imposes a small-rotationally-symmetric dimple very close to the tip of each initial surface. The numerical simulations in~\cite{GIKW} show that for Near Class initial data, these dimples disappear
as MCF progresses, and Type-II behaviors are found to occur. For solutions originating from Near Class perturbations, as from unperturbed initial data, dilations near the tip approach ``bowl solitons'' (see Figure~\ref{fig2}).

\begin{figure}[H]
	\includegraphics[width=0.65\textwidth]{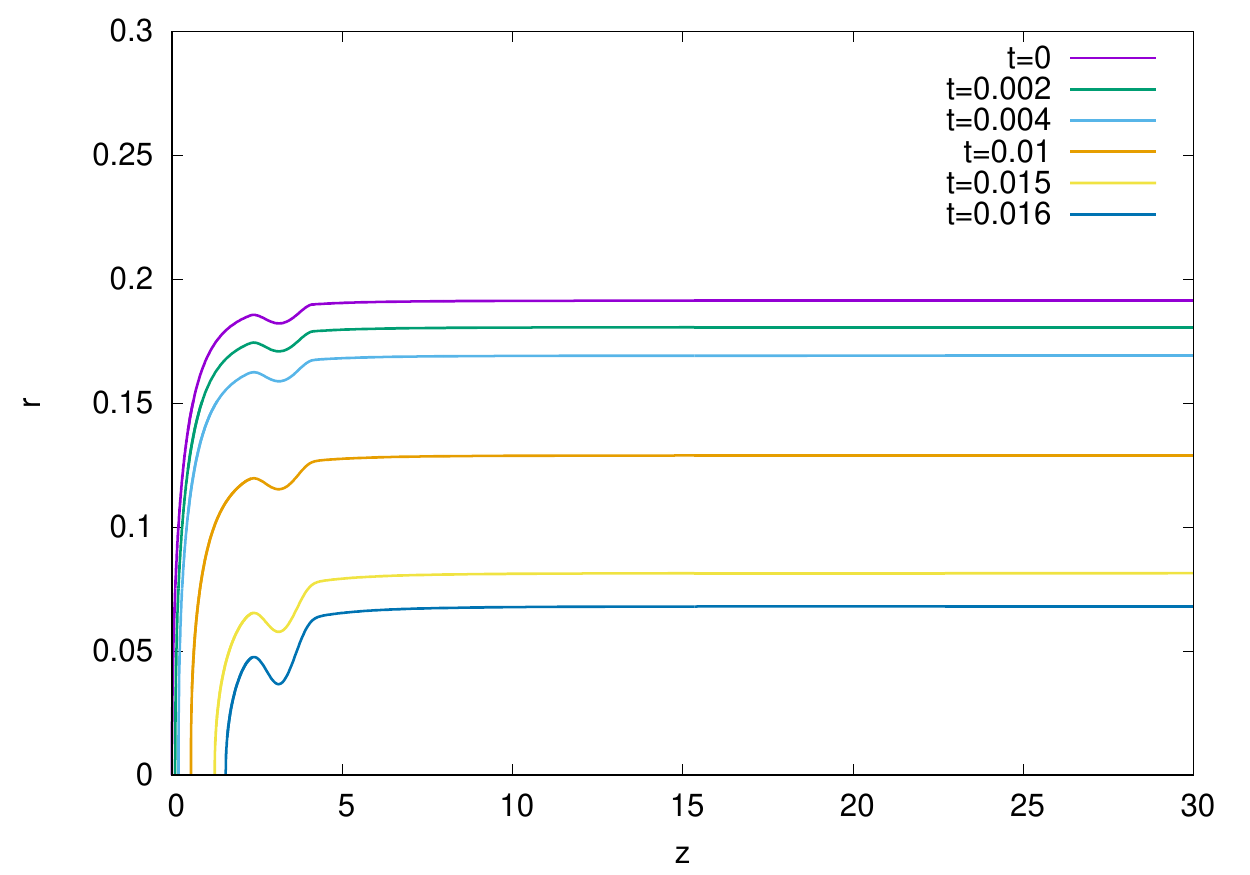}
	\caption{Numerical simulation of a MCF solution in the Far Class.
	Rotating each colored curve around the $z$-axis generates the hypersurface at the corresponding time.
	As the Type-I singularity develops, the dimple becomes a neckpinch.}			\label{fig3}
\end{figure}

The ``Far Class'' imposes a small rotationally-symmetric dimple in each initial surface relatively far from its tip. For this class of  solutions, MCF evolutions develop ``neckpinch'' singularities that exhibit Type-I
behaviors at the neck, with the curvatures at the tip remaining bounded (see Figure~\ref{fig3}).

As noted above, the main goal of this paper is to determine whether or not the results observed in the numerical simulations of~\cite{GIKW} remain true if we remove rotational symmetry. We have chosen to explore this issue by carrying out numerical simulations of MCF originating from initial data much like that of~\cite{GIKW}, but with sinusoidal angular dependence initially imposed. As we show below in Section~\ref{num_results}, the asymptotic behavior of these solutions is qualitatively unaffected  by this angular dependence. More precisely, the angular dependence disappears as time approaches the first singular time, and the singular behaviors approach those observed in either the rotationally-symmetric Near Class or Far Class, independent of the initially angular dependence.

While the need to work in a collection of independent but intersecting coordinate patches is very familiar in geometric analysis on manifolds, proceeding in this way is less familiar in numerical simulations of geometric evolutions. As discussed in \cite{GIKW}, 
the numerical simulations of MCF in that work requires the use of two coordinate patches.
\label{page:define-Coordinates} Using the coordinate $z$ for the direction parallel to the cylinder enveloping the embedded initial hypersurfaces, and the coordinate $r$ for the radial direction orthogonal to $z$ in the patch relatively close to the tip, the MCF evolutions are carried out in terms of the function $z(r,t)$. In the patch further away from the tip, this evolution is carried out in terms of the function $r(z,t)$.
Our numerical simulations in \cite{GIKW} function with these two patches overlapping, and have required the development of novel  numerical techniques that can handle this overlap.

In this work, we again focus on noncompact surfaces embedded in $\mathbb{R}^3$ and again work in two overlapping patches, here using $\theta$ as the angular coordinate. In the patch near the tip, a function $z(r,\theta,t)$ represents the MCF evolution, 
while in the patch farther from the tip, we use the function $r(z,\theta,t)$ instead.
The numerical techniques that we use here to handle MCF evolutions in the two overlapping patches are adapted from those in \cite{GIKW} to allow for this angular dependence. (A numerical analysis of non-rotational perturbations of higher-dimensional
embedded hypersurfaces will appear in a forthcoming work.)

In \cite{GIKW}, we have speculated that for rotationally-symmetric initial data at the transition between the  Near Class and the Far Class, MCF is likely to exhibit critical behavior, with Type-II degenerate-neckpinch singularity behavior occurring at the tip, and with curvature blowup rates consistent with such behavior. The asymptotic disappearance of angular dependence observed in the numerical simulations that we discuss in  this paper suggests that the existence of asymptotically rotationally-symmetric critical behaviors of MCF solutions is likely to occur for such non-rotationally-symmetric transitional initial embeddings as well.

Because the numerical simulations carried out here (as well as in~\cite{GIKW}) depend very much on the rigorous analysis of unperturbed MCF discussed in [\cite{IW} and~\cite{IWZ1}, we briefly review some results of those papers in Section~\ref{setup} and there set up the equations for the non-rotationally-symmetric case. In Section~\ref{num_method}, we write down the evolution equations for MCF of such non-rotationally-symmetric surfaces and discuss the computational techniques we use to carry out the numerical  simulations done here. In Section~\ref{num_results}, we present the results of these simulations, both for the Near Class and  the Far Class of MCF originating from non-rotationally-symmetric initial data. The conclusions are summarised in Section \ref{conclusion}.


\section{MCF of unperturbed initial data}\label{setup}

We briefly review some background analysis useful in studying the MCF solutions that are shown in~\cite{IW, IWZ1} to develop Type-II singularities. These flows evolve from what we call our ``unperturbed initial data'' --- complete rotationally-symmetric ``undimpled'' surfaces as shown in Figure~\ref{fig1}.

For any point $(x_0,~x_1,~x_2)\in\mathbb{R}^{3}$, we write $z = x_0$ and $r = \sqrt{x_1^2 + x_2^2}$. Let $\Gamma$ be a noncompact surface obtained by rotating the graph of $r(z)$, $a\leq z<\infty$, around the $z$-axis. The function $r(z)$ is assumed to be strictly concave, so that $\Gamma$ is strictly convex and $r$ is strictly increasing with $r(a)=0$ and $\lim\limits_{z\nearrow \infty}r(z) = r_0$, where $r_0$ is the radius of the enveloping cylinder. The function $r$ is assumed to be smooth except at $z=a$. The non-smoothness of $r$ is due to the particular choice of coordinates; if the time-dependent profile function $r(z,t)$ were inverted in a suitable way, then this irregularity would be removed. We call the point where $r=0$ the tip of the hypersurface.

Consider the class $\mathscr{G}$ of complete surfaces that are rotationally symmetric, (strictly) convex,\footnote{Throughout this paper, ``convex'' means ``strictly convex''.}
smooth graphs over a ball and are asymptotic to a cylinder. Then a MCF solution staring from any surface in this class remains embedded under the flow and moves towards its open end while remaining asymptotic to a shrinking cylinder \cite{SS14}. The solution disappears at spatial infinity at the same finite time $T$, called the vanishing time, that the cylinder collapses. This phenomenon can be viewed as a MCF singularity forming in finite time at spatial infinity. In fact, as shown in \cite{IW, IWZ1}, this singularity is Type-II: for each real number $\gamma\geq 1/2$, there exist MCF
solutions in the class $\mathscr{G}$ with the Type-II curvature blowup rate $(T-t)^{-(\gamma+1/2)}$ forming at the tip. For these solutions, the geometry near the tip is modelled by a bowl soliton, and the geometry near spatial infinity approaches a cylinder. Our analysis proceeds in two cases: $\gamma>\frac12$ and $\gamma=\tfrac12$.

If $\gamma>1/2$, we recall from \cite{IW} that it is useful to implement rescaled time and space variables defined by
\begin{align*}
	\tau = -\log(T-t), \qquad y = z(T-t)^{\gamma-1/2}, \qquad \phi(y,\tau) = r(x,t)(T-t)^{-1/2}.
\end{align*}
Let $\lambda:=-1/y$. Then the MCF equation of a rotationally symmetric surface in terms of $\lambda(\phi,\tau)$ is
\begin{align}\label{lambda}
	\left. \p_\tau \right\vert_\phi \lambda & = \frac{\lambda_{\phi\phi}-2\lambda^2_\phi/\lambda}{1+e^{2\gamma\tau}\lambda^2_\phi/\lambda^4} + \left(\frac{1}{\phi} - \frac{\phi}{2}\right)\lambda_\phi + \left(\gamma-\frac{1}{2}\right)\lambda,
\end{align}
where the notation $\left. \p_\tau \right\vert_\phi $ means taking the partial $\tau$-derivative while keeping $\phi$ fixed, $\phi\in(-\sqrt{2},\sqrt{2})$, and $\tau\geq\tau_0:=-\log(T-t_0)$ for an initial time $t_0$. Taking any real number $c>0$ 
(\emph{e.g.,} $c=1$) and letting $A:= c2^{\gamma-1/2}$, we define global initial data for \eqref{lambda} by smoothing the piecewise-smooth profile function
\begin{align}\label{id}
	\hat\lambda_0(\phi) := \left\{
	\begin{array}{lr}
		\begin{array}{l} -A + e^{-2\gamma\tau_0}F(\zeta)- e^{-2\gamma\tau_0}F(R_1) \\
			\quad + \left[ A -c\left(2-(R_1e^{-\gamma\tau_0})^2\right)^{\gamma-1/2} \right] 
		\end{array}, & 0\leq|\zeta|\leq R_1,\\ \\
		-c(2-\phi^2)^{\gamma-1/2}, & R_1e^{-\gamma\tau_0}\leq|\phi|<\sqrt{2},
	\end{array}
	\right.
\end{align}
where $\zeta=\phi e^{\gamma\tau}$, $R_1$ is some large constant, and $F$ is the unique solution of the \textsc{ode} initial value problem
\begin{align*}
	\frac{F_{\zeta\zeta}}{1+F^2_{\zeta}/A^4} + F_{\zeta}/\zeta = (\gamma-1/2)A,\quad\quad F(0)=F_{\zeta}(0)=0.
\end{align*}

If $\gamma=1/2$, we recall from \cite{IWZ1} that it is useful to implement rescaled variables defined by
\begin{align*}
	\tau  &= -\log(T-t),\qquad y = z + a \log(T-t),\qquad \phi(y,\tau) = r(x,t)(T-t)^{-1/2}.
\end{align*}
Let $\lambda:=-1/y$. Then the MCF equation of a rotationally symmetric surface in terms of $\lambda(\phi,\tau)$ is
\begin{align}\label{lambda-crit}
	\left. \p_\tau \right\vert_\phi \lambda & = \frac{\lambda_{\phi\phi}-2\lambda^2_\phi/\lambda}{1+e^{\tau}\lambda^2_\phi/\lambda^4} + \left(\frac{1}{\phi} - \frac{\phi}{2}\right)\lambda_\phi - a \lambda^2,
\end{align}
where $\phi\in(-\sqrt{2},\sqrt{2})$, and $\tau\geq\tau_0:=-\log(T-t_0)$ for an initial time $t_0$. Taking any real number $c>0$
(\emph{e.g.,} $c=1$) and letting $A:= c 2^{\gamma-1/2}$, we define global initial data for~\eqref{lambda-crit} by smoothing the piecewise-smooth profile function
\begin{align}\label{id-crit}
	\widehat\lambda_0(\phi) := \left\{
	\begin{array}{lr}
		-A + e^{-\tau_0}F(\zeta)- e^{-\tau_0}F(R_1)\\ \quad +A-\left(c-a\log(2-R_1^2e^{-\tau_0})\right)^{-1}, & 0\leqslant |\zeta|\leqslant R_1, \\ \\
		-1/\left(c-a\log(2-\phi^2)\right), & R_1e^{-\tau_0/2}\leqslant |\phi|<\sqrt{2},
	\end{array}
	\right.
\end{align}
where $\zeta=\phi e^{\gamma\tau}$, $R_1$ is some large constant, and $F$ is the unique solution of the \textsc{ode} initial value problem
\begin{align*}
	\frac{F_{\zeta\zeta}}{1+F^2_{\zeta}/A^4} + F_{\zeta}/\zeta = aA^2,\quad\quad F(0)=F_{\zeta}(0)=0.
\end{align*}

In both cases, the unperturbed initial surfaces are determined by the constants $(n, c, R_1, \tau_0)$. The initial data profile function $\hat\lambda_0$ 
defined in either \eqref{id} or \eqref{id-crit} is piecewise-smooth and continuous, which is sufficient for the numerical simulations based on the finite-differencing method that we use in this paper. In Section~\ref{num_method}, we explain how to specify the perturbed initial data used in our simulations.


\section{Numerical method}\label{num_method}

Our numerical method essentially consists of writing the mean curvature flow PDE in a standard parabolic form and then using a standard finite-difference method for parabolic equations.  However, for reasons to be explained below, we find it convenient to write the equation in two different ways and to use both of them in implementing our numerical simulations.

In cylindrical coordinates, the surface is given by specifying $r$ as a function of $z$ and $\theta$.  Under mean curvature flow, the function $r(t,z,\theta)$ satisfies the following evolution equation:
\begin{equation}
{\partial _t} r = {\frac {(1+ {r_z ^2}){r_{\theta \theta}} + ({r^2} + {r_\theta ^2}){r_{zz}} - 2 {r_\theta}{r_z}{r_{\theta z}} - {r^{-1}} {r_\theta ^2}} {{r_\theta ^2} + {r^2}(1 + {r_z ^2})}} \; - \; {\frac 1 r} \; \; \; .
\label{mcfzth}
\end{equation}
Equation~(\ref{mcfzth}) becomes singular where $r=0$.  In our previous paper treating rotationally-symmetric surfaces, we have dealt with this issue by treating a part of the surface near $r=0$ by instead specifying $z$ as a function of $r$. The generalization of this technique to the case of non-rotationally symmetric surfaces would be to specify $z$ as a function of $r$ and $\theta$.  However, that would not be suitable, because the $(r,\theta)$ coordinate system becomes singular at $r=0$.  Instead, in a neighborhood of the tip, we use a Cartesian coordinate system, and specify the surface by writing $z$ as a function of $x$ and $y$.  Under mean curvature flow, the function $z(t,x,y)$ satisfies the following evolution equation:
\begin{equation}
{\partial _t} z = {\frac {(1+{z_y ^2}){z_{xx}} + (1+{z_x ^2}){z_{yy}} - 2 {z_x}{z_y}{z_{xy}}} {1 +{z_x ^2} + {z_y ^2}}}.
\label{mcfxy}
\end{equation}

The numerical evolution is done using Euler's method: that is, given the surface at $t$,
we use equation~\eqref{mcfzth} to evolve $r$ to time $t+\Delta t$ by $r(t+\Delta t) = r(t) + \Delta t {\partial _t} r$.  Here 
${\partial _t} r$ represents the right-hand side of equation~\eqref{mcfzth} evaluated at time $t$ by using the finite-difference approximations for spatial derivatives given below. Correspondingly, Euler's method is applied to equation~\eqref{mcfxy} to evolve $z$ from time $t$ to time $t+\Delta t$.

Our finite-difference approximations are standard second-order centered differences.  That is, for the function $r(z,\theta)$, we let $r^i _k$ denote $r(i \Delta z, k \Delta \theta)$.  Then our finite-difference approximations are
\begin{eqnarray}
{r_z} &=& {\frac {{r^{i+1} _k} - {r^{i-1} _k}} {2 \Delta z}},
\label{fd1a}
\\
{r_\theta} &=& {\frac {{r^i _{k+1}} - {r^i _{k-1}}} {2 \Delta \theta}},
\label{fd1b}
\\
{r_{zz}} &=& {\frac {{r^{i+1} _k} + {r^{i-1} _k} - 2 {r^i _k}} {{\Delta z}^2}},
\label{fd1c}
\\
{r_{\theta \theta}} &=& {\frac {{r^i _{k+1}} + {r^i _{k-1}} - 2{r^i _k}} {{\Delta \theta }^2}},
\label{fd1d}
\\
{r_{z\theta}} &=& {\frac {({r^{i+1} _{k+1}} + {r^{i-1} _{k-1}}) - ({r^{i+1} _{k-1}} + {r^{i-1} _{k+1}} ) } {4 \Delta z \Delta \theta}}.
\label{fd1e}
\end{eqnarray}

Similarly, for the function $z(x,y)$, we let $z^i _k$ denote $z(i \Delta x, k \Delta y)$.  Then our finite-difference approximations are
\begin{eqnarray}
{z_x} &=& {\frac {{z^{i+1} _k} - {z^{i-1} _k}} {2 \Delta x}},
\label{fd2a}
\\
{z_y} &=& {\frac {{z^i _{k+1}} - {z^i _{k-1}}} {2 \Delta y}},
\label{fd2b}
\\
{z_{xx}} &=& {\frac {{z^{i+1} _k} + {z^{i-1} _k} - 2 {z^i _k}} {{\Delta x}^2}},
\label{fd2c}
\\
{z_{yy}} &=& {\frac {{z^i _{k+1}} + {z^i _{k-1}} - 2{z^i _k}} {{\Delta y }^2}},
\label{fd2d}
\\
{z_{xy}} &=& {\frac {({z^{i+1} _{k+1}} + {z^{i-1} _{k-1}}) - ({z^{i+1} _{k-1}} + {z^{i-1} _{k+1}} ) } {4 \Delta x \Delta y}}.
\label{fd2e}
\end{eqnarray}

Euler's method is slow for parabolic equations, because numerical stability requires that, in Cartesian coordinates, the time step $\Delta t$ be smaller than a number of order 1 times the smaller of the squares of the space steps ${\Delta x}^2$ and ${\Delta y}^2$.  In $(z,\theta)$ coordinates, the bounds on permissible $\Delta t$ are given by ${\Delta z}^2$ and ${(r \Delta \theta )}^2$.  This makes the evolution quite slow when we approach a neck pinch, since $r \to 0$ there.  If we were to use $(r,\theta)$ coordinates rather than $(x,y)$ coordinates, the bounds on permissible $\Delta t$ would be ${\Delta r}^2$ and ${(r\Delta \theta )}^2$.  Since $r \to 0$ at the origin in polar coordinates, this would give rise to a very small $\Delta t$.  This is another way of saying why Cartesian coordinates are preferable to polar coordinates in this case.

Our overlap method requires that we use two grids, each of which is a manifold with boundary.  The boundary points of one grid correspond to interior points of the other grid, and their values are set using interpolation in the other grid.  In our previous paper, one grid had a single coordinate $r$ and a single boundary point given by $r_{\rm max}$ while the other grid had a single coordinate $z$ with a single boundary point $z_{\rm min}$.  Here the boundary of our $(z,\theta)$ grid is the circle $(z_{\rm min},\theta)$, while the boundary of our $(x,y)$ grid is a square of side length $L$ with center at the origin.  As with our previous method, as the surface evolves, we adjust the parameters (in this case $z_{\rm min}$ and $L$) so that each boundary remains within the interior of the other grid.

We end this section by explaining how we choose the perturbed initial data sets for the numerical simulations. We first recall that for an unperturbed solution constructed in \cite{IW} (or \cite{IWZ1}), respectively, its initial data set  is obtained by joining a scaled bowl soliton to a cylinder at spatial infinity, and is defined precisely in the rescaled coordinates by equation \eqref{id} (or \eqref{id-crit}, respectively). Expressed in terms of the unscaled $(z,r)$-coordinates, the ODE for the bowl soliton takes the form
\begin{equation}\label{bowl-ode}
z_{rr} = (1+z_r^2)\left(\beta - \frac{z_r}{r}\right).
\end{equation}

For each choice of $\gamma$, we fix $c=1$ and choose $\tau_0$ large (in all the numerical simulations $\gamma =3/4$ and $\tau_0=4$) so that the matched asymptotics explained in Section \ref{setup} make sense. We then choose $\beta$ according to $\beta=c^{-1}(\gamma-1/2)2^{-(\gamma-1/2)} e^{-(\gamma+1/2)\tau_0}$, and numerically integrate equation \eqref{bowl-ode} outward from $r=0$ to value $r_1$, which is obtained by writing $R_1$ (chosen to be $e^{\gamma \tau_0/2}$ in our numerical simulations) in the unscaled $r$-coordinate. For all $r>r_1$, we use the following analytic formula which follows from rewriting the second equation in \eqref{id} in the $(r,z)$-coordinates:
\begin{equation}\label{formula}
z(r) = z(r_1)+\frac{1}{c}\left[ \left( 2e^{-\tau_0}-r^2\right) ^{\frac{1}{2}-\gamma} - \left(2e^{-\tau_0}-r_1^2\right)^{\frac{1}{2}-\gamma} \right].
\end{equation}

We choose the patch where we write $z$ as a function of $(x,y)$ to be a square whose sides have length $10 {r_1}$.  In this patch, we have $r={\sqrt {{x^2}+{y^2}}}$.  For $r < {r_1} $, we specify $z$ by the results of the numerical integration of \eqref{bowl-ode}.  For $r \ge {r_1}$, $z$ is given by formula \eqref{formula}.

In the patch where $r$ is written as a function of $(z,\theta)$, we use the formula obtained by inverting expression \eqref{formula}.

For our ``Near Class'' of distorted surfaces, we change the surface only in the $(x,y)$ patch. We choose a number $r_m$ and distort the surface only where $r < {r_m}$.  In this region, we choose two amplitudes $a_0$ and $a_1$ and change $z$ by 
\begin{equation}
z \mapsto z + {a_0} (1 + {a_1} x y ) ({r^2} -  {r_m ^2})
\end{equation}

For our ``Far Class'' of distorted surfaces, we change the surface only in the $(z,\theta)$ patch. We choose an amplitude $a_0$, and two values of $z$: $z_a$ and $z_b$. We distort the surface only on the interval $(z_a, z_b)$.  On this interval, we define the quantity $F$ by 
\begin{equation}
F \equiv \left ( 1 + {\textstyle {\frac 1 4}} \cos (n \theta ) \right ) \sin \left ( {\frac {\pi (z - {z_b})} {{z_a}-{z_b}}} \right ) \; \; \; .
\end{equation}
We then distort the surface by 
\begin{equation}
r \mapsto r (1-{a_0}{F^2}).
\end{equation} 
Here $n$ is an integer, which we choose to be even so that the neck pinch will occur at $r=0$.  An odd $n$ would
give rise to an off-center neck pinch, which results in a coordinate singularity when the shrinking surface intersects $r=0$.

\section{Numerical results}\label{num_results}

\begin{figure}[H]
	\includegraphics[width=0.65\textwidth]{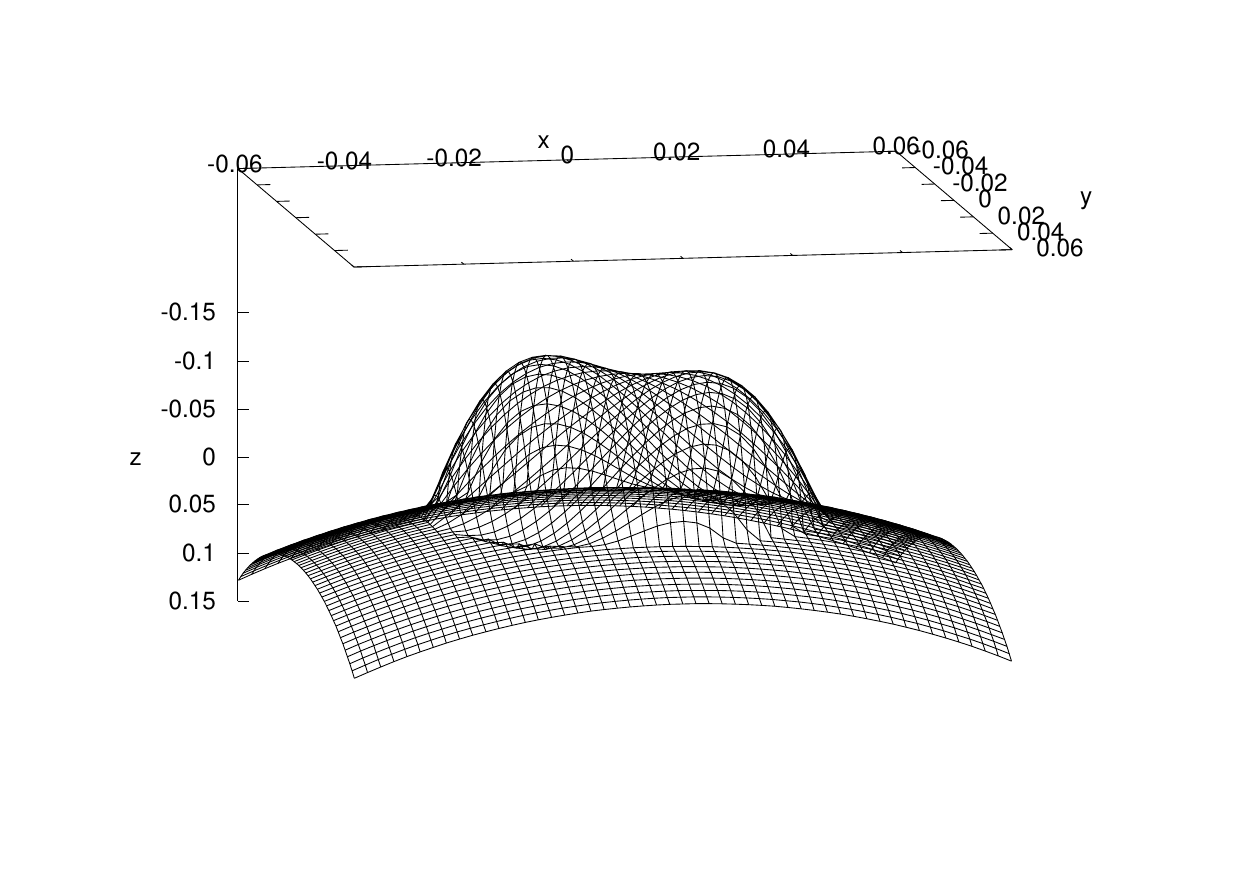}
	\caption{A simulated three-dimensional graph of angular dependent initial embedding for Near Class data in a neighborhood of the tip.}				\label{initcoarse}
\end{figure}

\begin{figure}[H]
	\includegraphics[width=0.65\textwidth]{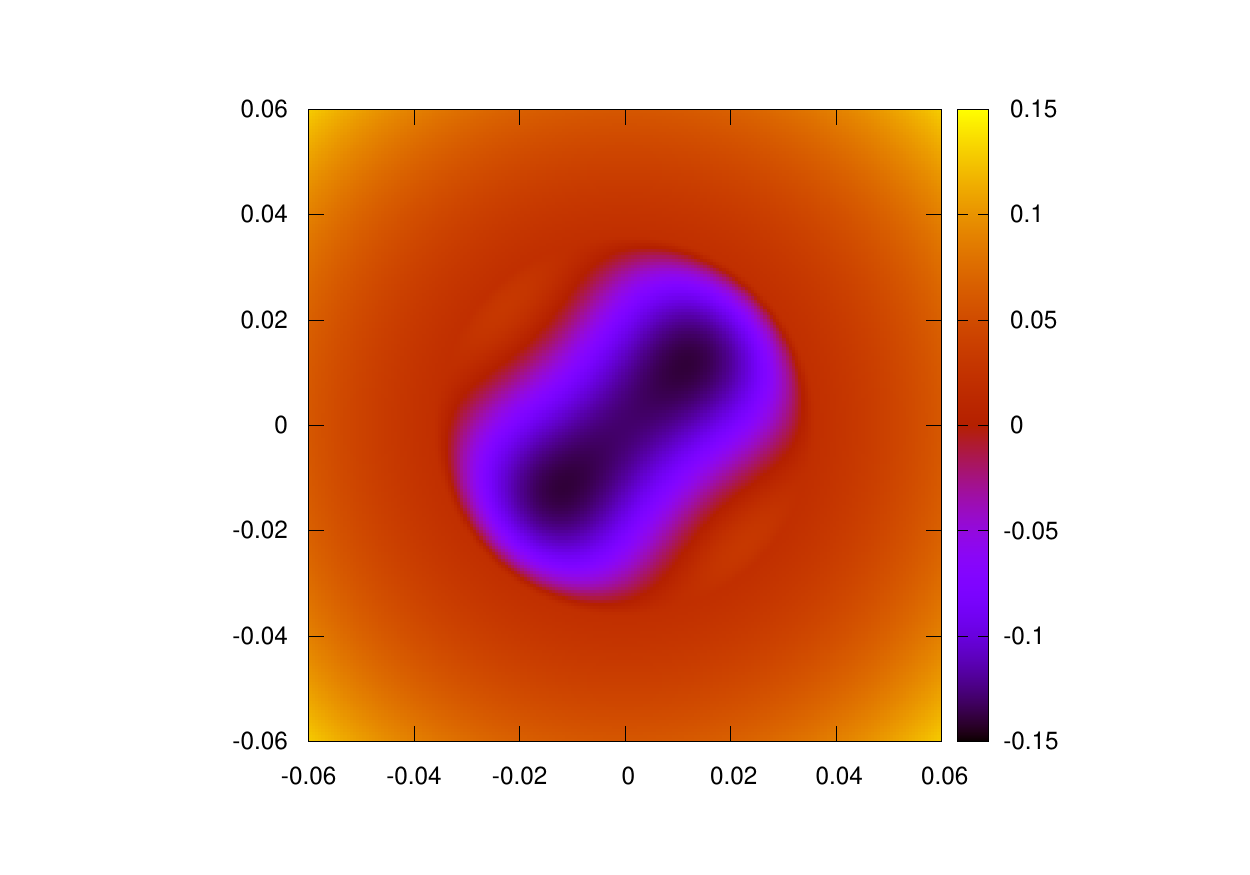}
	\caption{A colored ``heat map'' of the initial data represented in Figure~\ref{initcoarse}. The colors correlate with the graph height.}				\label{init}
\end{figure}

\begin{figure}[H]
	\includegraphics[width=0.65\textwidth]{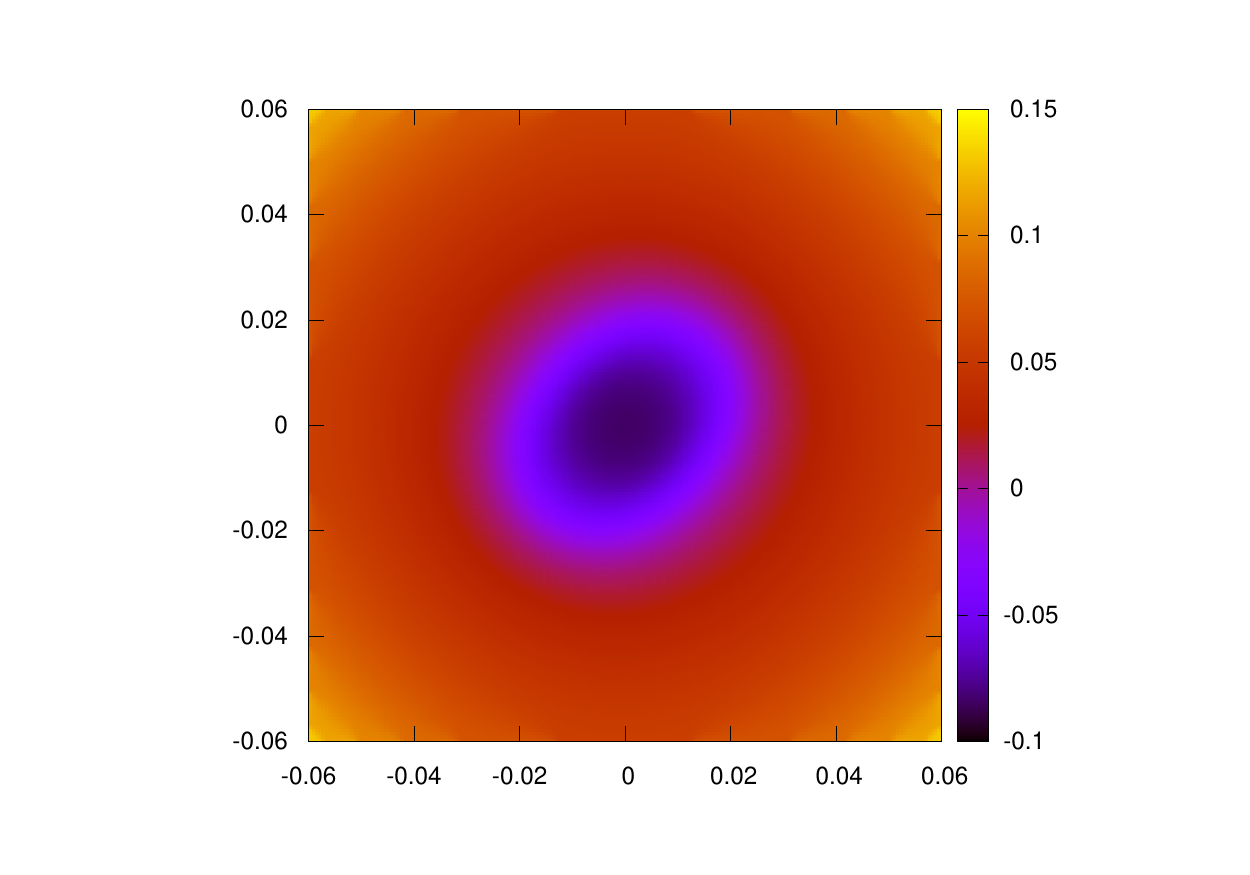}
	\caption{A colored heat map of the mean curvature flow of the initial data depicted in Figure~\ref{init}. It is clearly becoming rounder.}				\label{t1}
\end{figure}

\begin{figure}[H]
	\includegraphics[width=0.65\textwidth]{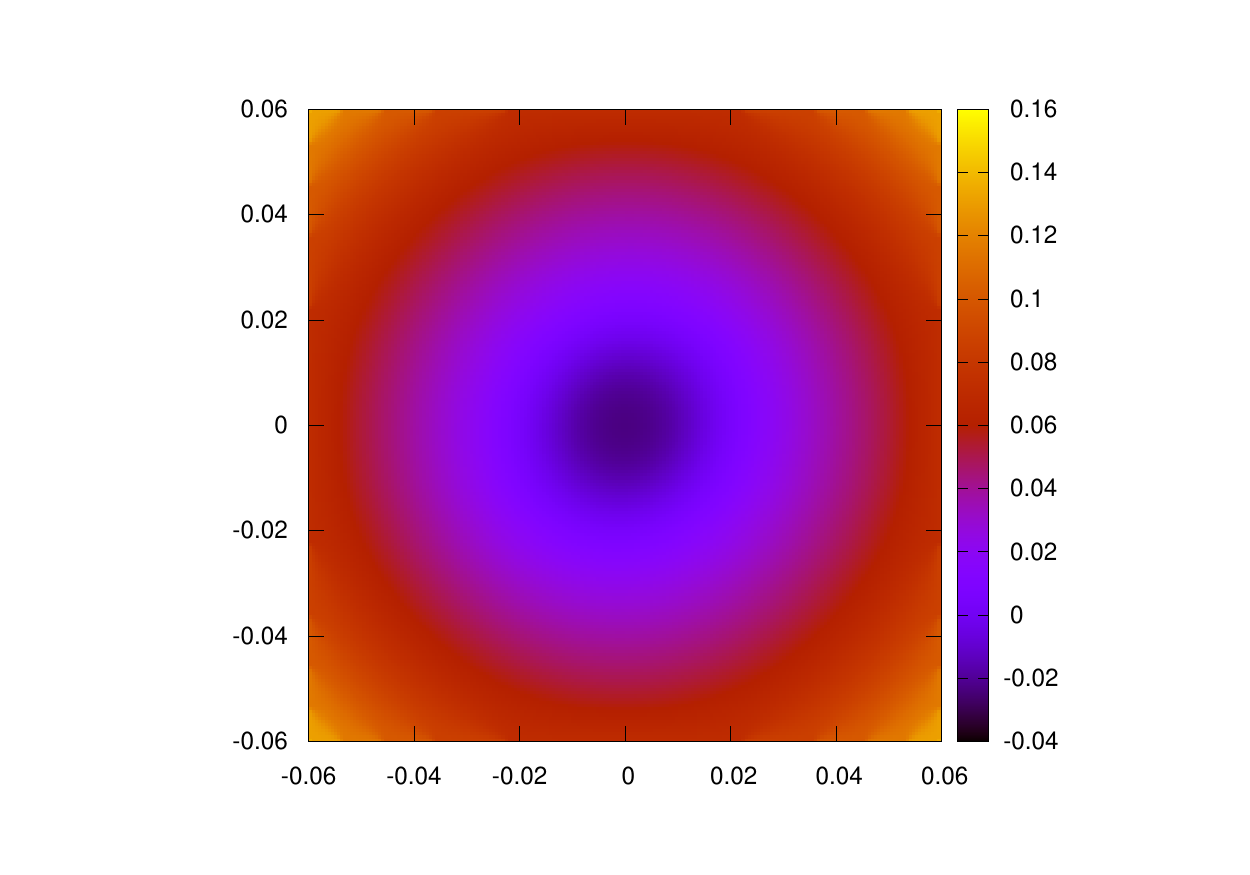}
	\caption{A colored heat map of the mean curvature flow with initial data from Figure~\ref{init} at a later time. This time, the angular dependence is gone.}				\label{t2}
\end{figure}

We present our results primarily graphically. In the first pair of graphs – Figure~\ref{initcoarse} and Figure~\ref{init} – we represent the initial hypersurface in a neighborhood of the tip. The angular dependent dimple is represented as a simulated three-dimensional graph in Figure~\ref{initcoarse} and as a color coded ``heat map'' in Figure~\ref{init}. The ``rounding'' of the mean curvature flow evolving from the initial hypersurface of Figure~\ref{init} is in a heat map at a later time in Figure~\ref{t1}, and then even later in Figure~\ref{t2}, when the hypersurface becomes almost entirely rotationally symmetric.

We do not also include graphs of the evolution of the curvature at the tip for angular dependent Near Class data, but we do note that such graphs (similar to sub-figures (B) and (C) in Figures 2--5 in \cite{GIKW}) do show that, for such mean curvature flows, the numerical evidence supports the contention that Type-II singularities occur at the tip for perturbed initial data, even with angular dependence.

\begin{figure}[H]
	\includegraphics[width=0.65\textwidth]{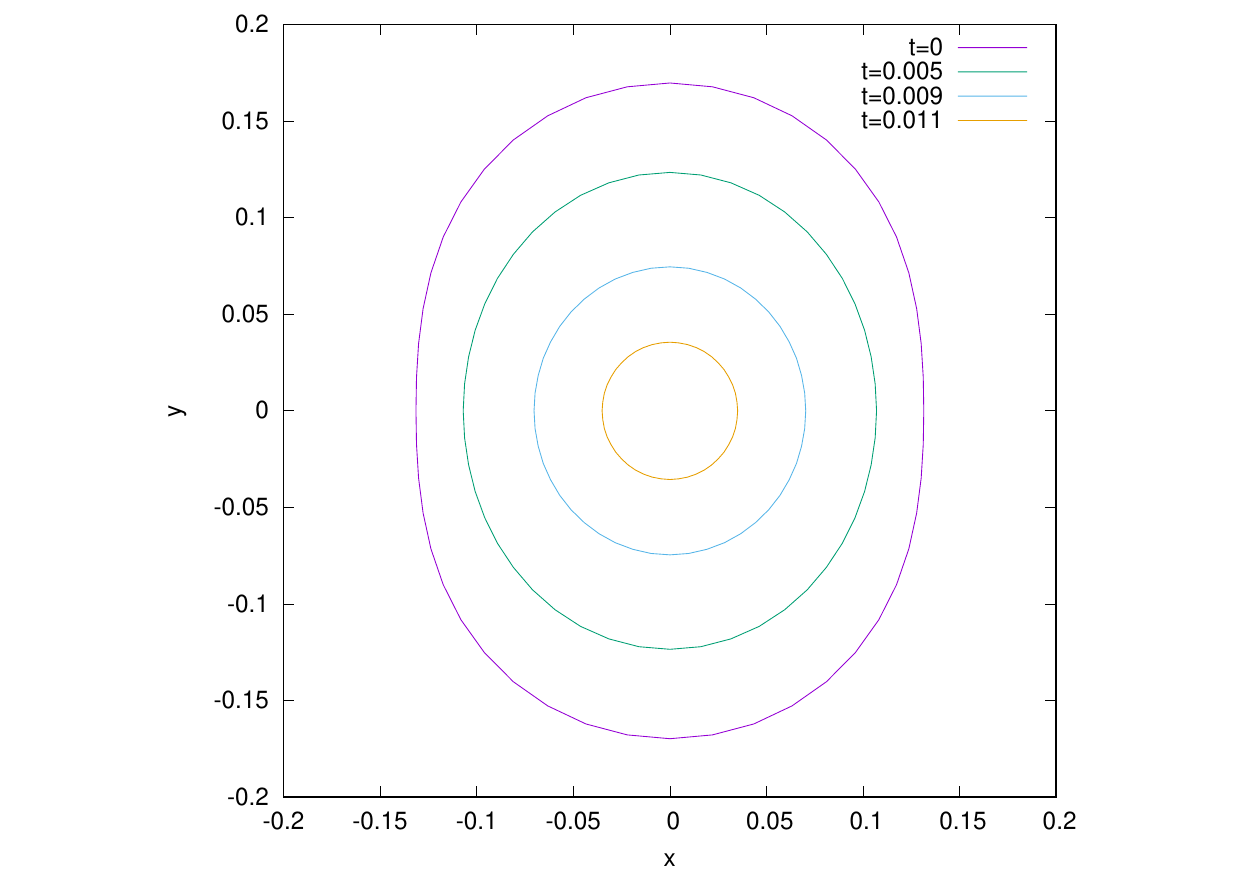}
	\caption{Vertical cross-sections of the mean curvature flow for far class initial data at successive times, with $\cos(2\theta)$ angular dependence. These cross-sections are chosen to be the location of the developing neck pinch. The loss of the angular dependence is evident.}				\label{neckcos2theta}
\end{figure}

\begin{figure}[H]
	\includegraphics[width=0.65\textwidth]{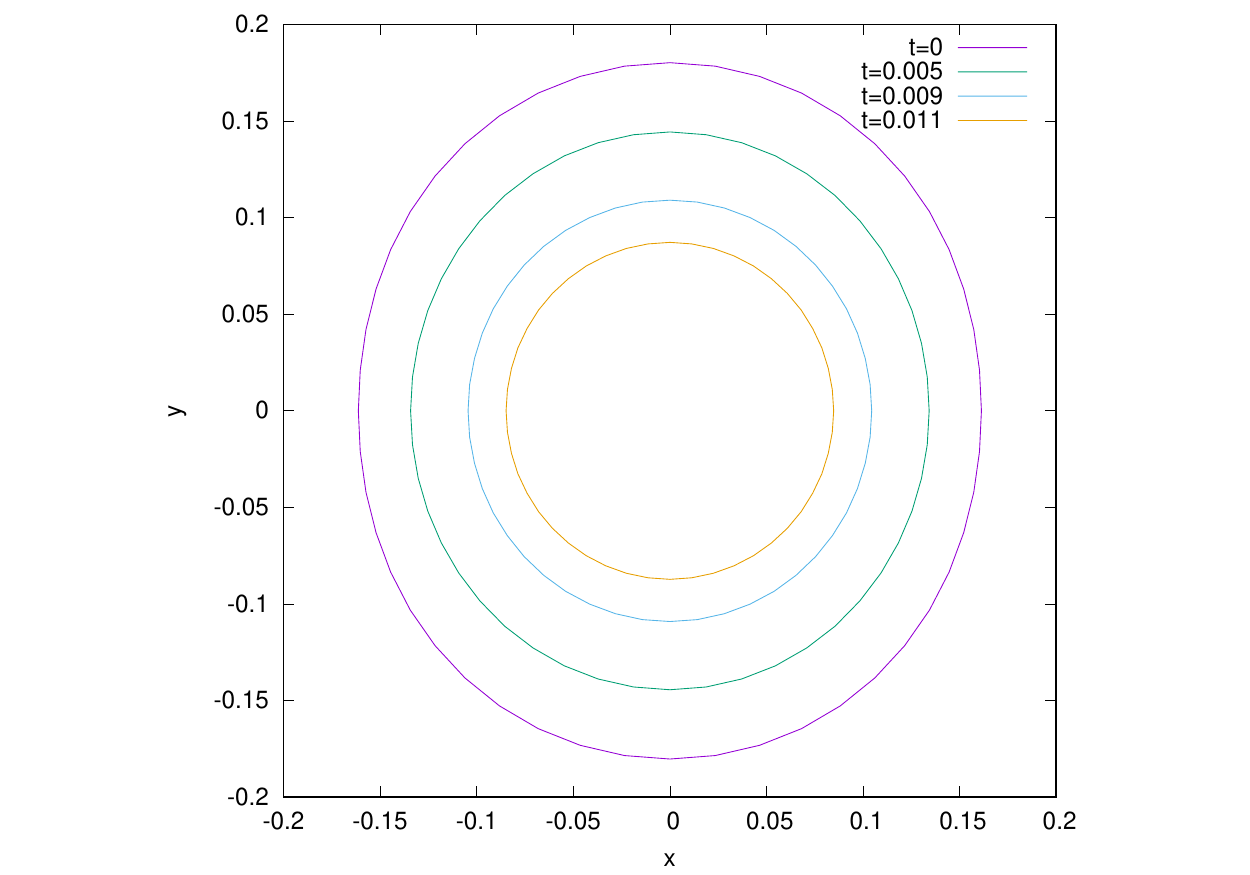}
	\caption{Vertical cross-sections of the MCF for the same initial data in Figure~\ref{neckcos2theta}, but with the cross-sections obtained to the left of the neck pinch.}				\label{leftofneckcos2theta}
\end{figure}

\begin{figure}[H]
	\includegraphics[width=0.65\textwidth]{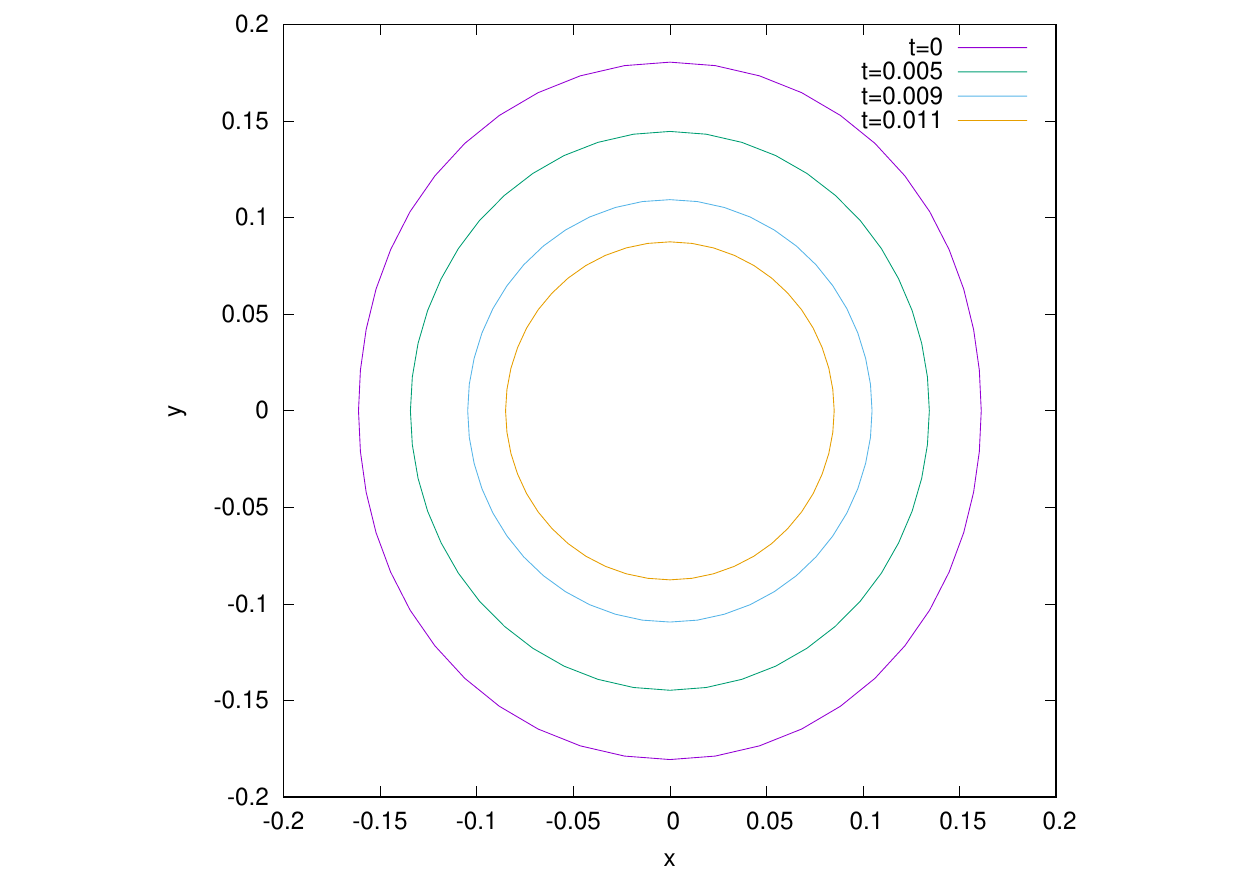}
	\caption{Vertical cross-sections of the MCF for the same initial data in Figure~\ref{neckcos2theta}, with the cross-sections obtained to the right of the neck pinch.}				\label{rightofneckcos2theta}
\end{figure}

In Figures~\ref{neckcos2theta} through~\ref{rightofneckcos2theta}, we graphically represent the mean curvature flow for initial data with $\cos(2\theta)$ angular dependence. In Figure~\ref{neckcos2theta}, each of the four colored curves represents a vertical cross-section of the mean curvature flow at the location of the evolving neck pinch for far data at four successive times. The loss of the angular dependence is evident. In Figure~\ref{leftofneckcos2theta}, a similar graphical representation of vertical cross-sections of the mean curvature flow is depicted a bit to the left of the neck pinch, while in Figure~\ref{rightofneckcos2theta}, the same is depicted for vertical cross-sections a bit to the right of the neck pinch. In all of these figures, the rounding is evident.

\begin{figure}[H]
	\includegraphics[width=0.65\textwidth]{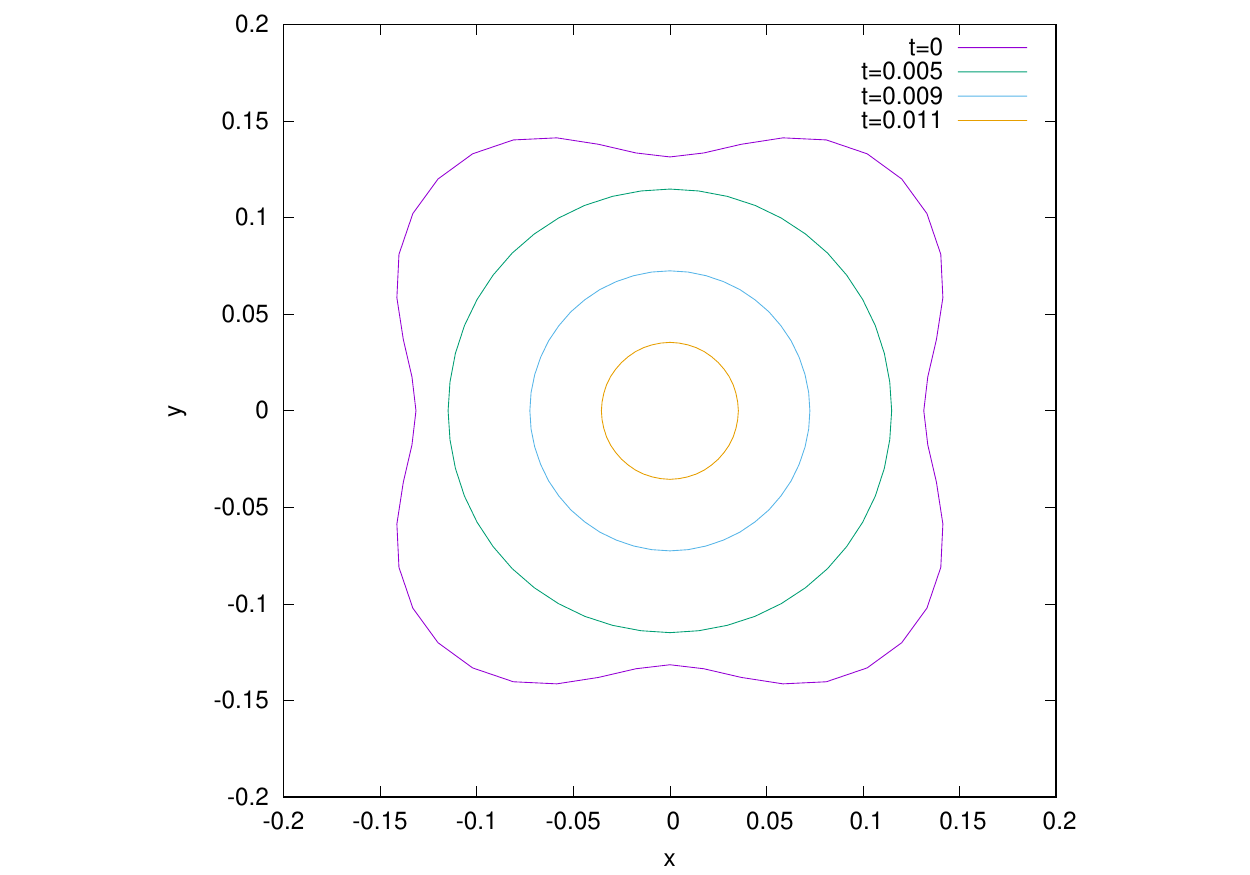}
	\caption{Vertical cross-sections similar to Figure~\ref{neckcos2theta}, but with imposed angular dependence of $\cos(4\theta)$.}				\label{neckcos4theta}
\end{figure}

\begin{figure}[H]
	\includegraphics[width=0.65\textwidth]{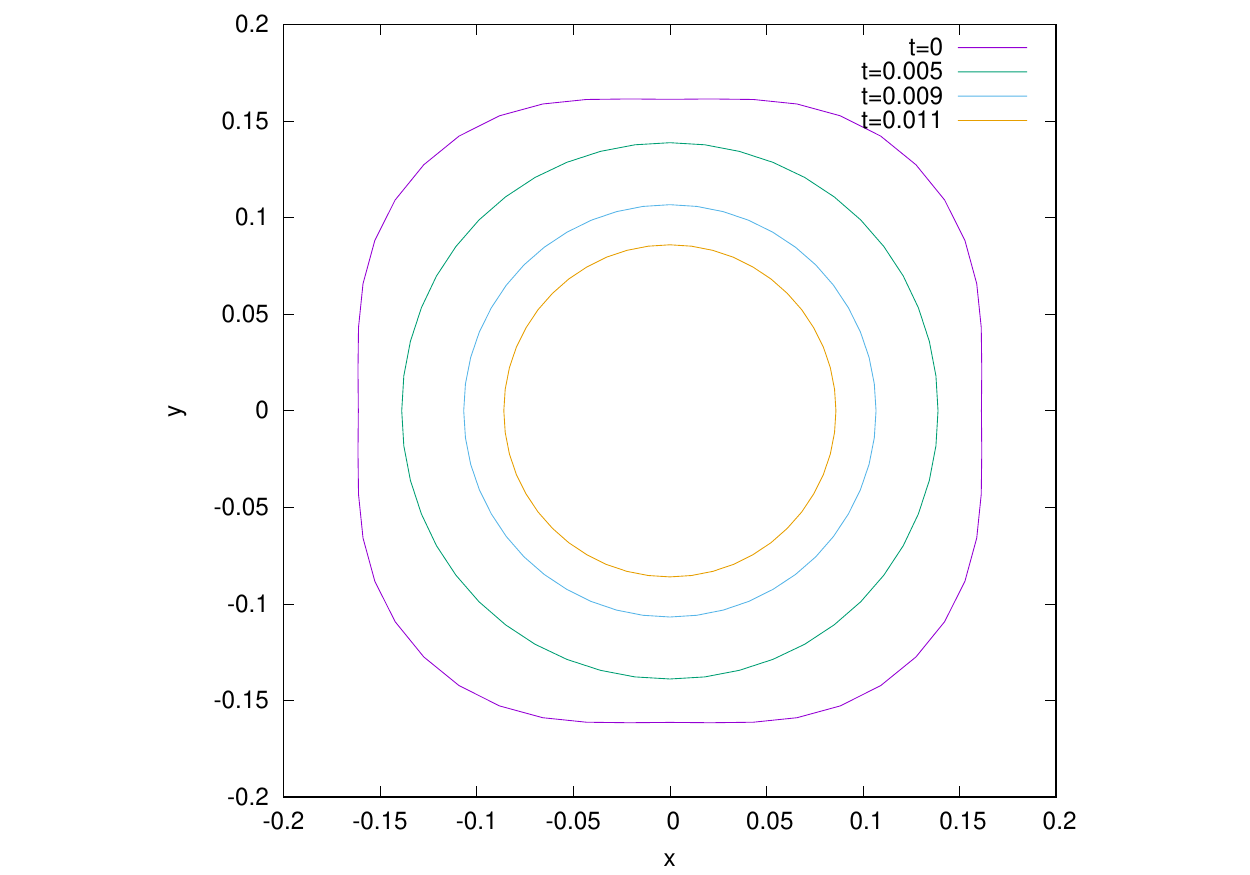}
	\caption{Vertical cross-sections similar to Figure~\ref{leftofneckcos2theta}, with imposed angular dependence $\cos(4\theta)$.}				\label{leftofneckcos4theta}
\end{figure}

\begin{figure}[H]
	\includegraphics[width=0.65\textwidth]{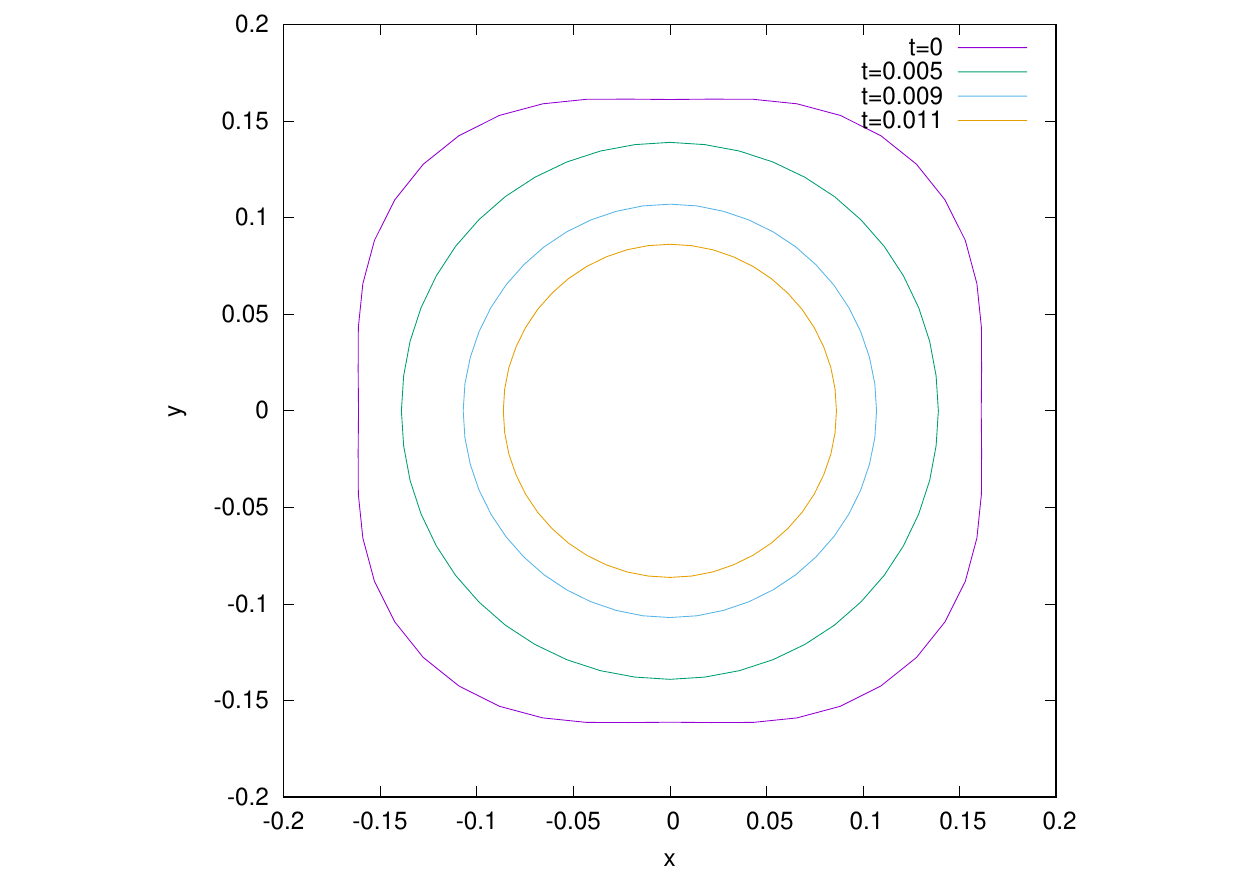}
	\caption{Vertical cross-sections similar to Figure~\ref{rightofneckcos2theta}, with imposed angular dependence $\cos(4\theta)$.}				\label{rightofneckcos4theta}
\end{figure}

Figures~\ref{neckcos4theta} through~\ref{rightofneckcos4theta} are very similar to Figures~\ref{neckcos2theta} through~\ref{rightofneckcos2theta}, in their depiction of vertical cross-sections of the mean curvature flow for far data, but in this case, the imposed angular dependence is $\cos(4\theta)$. Consequently, the initial data shows larger angular dependence. But again, we see that as the flow evolves, these cross-sections become increasingly round.

\section{Conclusions}\label{conclusion}

These numerical simulations indicate that mean curvature evolutions of angular dependent perturbations asymptotically approach rotationally symmetric solutions. They do not prove that the mean curvature flow for general perturbations of the initial embedded hypersurfaces considered in \cite{IW,IWZ1} always have the same Type-II singularity behavior at the tip as is found in \cite{IW,IWZ1}. However, the numerical simulations carried out here together with those in \cite{GIKW} provide strong evidence that this singularity behavior is indeed stable. Mathematical proof of such stability deserves exploration.


\section*{Acknowledgements}
DG thanks the National Science Foundation for support in PHY-1806219 and PHY-2102914. JI thanks the National Science Foundation for support in grant PHY-1707427. DK thanks the Simons Foundation for support in Award 635293. HW thanks the Australian Research Council for support in DE180101348.



\bibliography{mcf_numerical_2}

\begin{bibdiv}
\begin{biblist}

\bib{GIKW}{article}{
      author={Garfinkle, David},
      author={Isenberg, James},
      author={Knopf, Dan},
      author={Wu, Haotian},
       title={A numerical stability analysis of mean curvature flow of
  noncompact hypersurfaces with type-{II} curvature blowup},
        date={2021},
        ISSN={0951-7715},
     journal={Nonlinearity},
      volume={34},
      number={9},
       pages={6539\ndash 6560},
  url={https://doi-org.ezproxy.library.sydney.edu.au/10.1088/1361-6544/ac15a9},
}

\bib{IW}{article}{
      author={Isenberg, James},
      author={Wu, Haotian},
       title={Mean curvature flow of noncompact hypersurfaces with {T}ype-{II}
  curvature blow-up},
        date={2019},
        ISSN={0075-4102},
     journal={J. Reine Angew. Math.},
      volume={754},
       pages={225\ndash 251},
         url={https://doi.org/10.1515/crelle-2017-0019},
      review={\MR{4000574}},
}

\bib{IWZ1}{article}{
      author={Isenberg, James},
      author={Wu, Haotian},
      author={Zhang, Zhou},
       title={Mean curvature flow of noncompact hypersurfaces with type-{II}
  curvature blow-up. {II}},
        date={2020},
        ISSN={0001-8708},
     journal={Adv. Math.},
      volume={367},
       pages={107111, 44},
         url={https://doi.org/10.1016/j.aim.2020.107111},
      review={\MR{4078822}},
}

\bib{IWZ2}{article}{
      author={Isenberg, J.},
      author={Wu, H.},
      author={Zhang, Z.},
       title={On the precise asymptotics of {T}ype-{II}b solutions to mean
  curvature flow},
     journal={Trans. Amer. Math. Soc. Ser. B},
       pages={to appear},
}

\bib{SS14}{article}{
      author={S\'{a}ez, Mariel},
      author={Schn\"{u}rer, Oliver~C.},
       title={Mean curvature flow without singularities},
        date={2014},
        ISSN={0022-040X},
     journal={J. Differential Geom.},
      volume={97},
      number={3},
       pages={545\ndash 570},
         url={http://projecteuclid.org/euclid.jdg/1406033979},
      review={\MR{3263514}},
}

\end{biblist}
\end{bibdiv}

\end{document}